\documentclass[12pt,a4paper,leqno]{amsart}
\usepackage{amsfonts,amssymb,amsmath,amsthm,graphicx}
\usepackage{a4wide}
\usepackage{mathrsfs}
\usepackage{epic,eepic}

\newtheorem{theorem}{Theorem}[section]
\newtheorem{proposition}[theorem]{Proposition}
\newtheorem{corollary}[theorem]{Corollary}
\newtheorem{lemma}[theorem]{Lemma}
\newtheorem{definition}[theorem]{Definition}
\theoremstyle{remark}
\newtheorem{remark}[theorem]{Remark}
\numberwithin{equation}{section}
\newtheorem{problem}{Problem}

\newcommand{\R}{{\mathbb R}}

\renewcommand{\L}{\mathscr{L}}
\renewcommand{\H}{{\mathcal H}}
\newcommand{\eps}{{\varepsilon}}
\newcommand{\dist}{\mathrm{dist}}
\newcommand{\const}{\mathrm{const.}}
\newcommand{\supp}{\mathrm{supp}\:}
\newcommand{\ul}{\underline}
\newcommand{\ol}{\overline}

\newcommand{\tx}[1]{\mbox{\;{#1}\;}}
\newcommand{\beqn}{\begin{eqnarray}}
\newcommand{\eeqn}{\end{eqnarray}}   
\newcommand{\beq}{\begin{eqnarray*}}
\newcommand{\eeq}{\end{eqnarray*}}
\newcommand{\hsub}{\underline{h}}
\newcommand{\Hsub}{\underline{H}}
\newcommand{\hsup}{\overline{h}}
\newcommand{\Hsup}{\overline{H}}
\newcommand{\usub}{\underline{u}}
\newcommand{\usup}{\overline{u}}
\newcommand{\vsub}{\underline{v}}

\begin{document}

\title["Boundary blowup" type sub-solutions]
{"Boundary blowup" type sub-solutions to semilinear elliptic equations with Hardy potential}
\author{Catherine Bandle, Vitaly Moroz and Wolfgang Reichel}
\address{C. Bandle \hfill\break Mathematisches Institut, Universit\"at
  Basel \hfill\break
  Rheinsprung 21, CH-4051 Basel, Switzerland}
  \email{C.Bandle@gmx.ch}
\address{V. Moroz \hfill\break Department of Mathematics, University
  of Bristol \hfill\break University Walk, Bristol BS8 1TW, UK}
  \email{V.Moroz@bristol.ac.uk}
\address{W. Reichel \hfill\break Mathematisches Institut,
  Universit\"at Gie{\ss}en \hfill\break Arndtstr. 2, D-35392 Gie{\ss}en, Germany}
  \email{Wolfgang.Reichel@math.uni-giessen.de}

\date{\today}

\begin{abstract}Semilinear elliptic equations which give rise to solutions
  blowing up at the boundary are perturbed by a Hardy
  potential $\mu/\delta(x,\partial\Omega)^2$. The size of this potential effects the
  existence of a certain type of solutions (large solutions): if $\mu$
  is too small, then no large solution exists. The presence of the Hardy
  potential requires a new definition of large solutions, following
  the pattern of the associated linear problem. Nonexistence and
  existence results for different types of solutions will be given. Our
  considerations are based on a Phragmen-Lindel\"of type theorem which
  enables us to classify the solutions and sub-solutions according to their behavior near the
boundary. Nonexistence follows from this principle together with the
Keller-Osserman upper bound.
The existence proofs rely on sub- and super-solution techniques
and on estimates for the Hardy constant derived in Marcus, Mizel and Pinchover \cite{MMP}.
\end{abstract}


\maketitle

\section{Introduction}

On bounded smooth domains $\Omega\subset\R^N, N\geq 2$ we study the
existence and non-existence of positive solutions and
sub-solutions to semilinear elliptic equations of the form
\begin{equation}\label{*}
-\Delta u-\frac{\mu}{\delta^2}u+\frac{u^p}{\delta^s}=0\quad\mbox{in }\Omega,
\end{equation}
where $\mu,s\in \R$ and $p>1$ are given constants and
$$\delta(x):=\dist(x,\partial\Omega).$$
There are two competing ingredients in \eqref{*}, namely the nonlinear problem
\beq
{\bf \mbox{(N)}}\quad\quad\quad\quad -\Delta U+\frac{U^p}{\delta^s}=0\quad\mbox{in }\Omega
\eeq
and the linear problem
\beq
{\bf \mbox{(L)}}\quad\quad\quad\quad-\Delta h-\frac{\mu}{\delta^2}h=0 \tx{in} \Omega.
\eeq
The nonlinear problem (N) has received a lot of attention in recent years,
cf. \cite{MV} and the references cited therein.
For $s<2$ it possesses a  maximal solution which is larger than
any other solution in $\Omega$. This solution behaves like
$c_{p,s}\delta(x)^\frac{s-2}{p-1}$. Since it tends to $+\infty$  as
$x$ approaches the boundary, it became common to call such solutions
{\em boundary blow-up} solutions or simply {\em large} solutions. For
$s>2$ only the trivial solution exists.
It follows from the Keller-Osserman upper bound given in Section \ref{kob}.
A related nonexistence result is found in \cite{Veron}. There problem (N) is considered
in the unit ball of $\R^N$ with $N\ge 3$, $p=\frac{N+2}{N-2}$ and $s\ge 2$.

\medskip
The linear problem (L) has been studied recently in \cite{BrMa97} and
in \cite{MMP} in connection  with Hardy's inequality. In this paper we are
interested only in positive solutions of (L). We shall call them
 {\em harmonics}. The concept of {\em sub-} and {\em
  super-harmonics} is understood in the usual pointwise
sense. It makes sense to extend the concept of
(sub-/super-)harmonics to {\em local} (sub-/super-)harmonics, which are
 defined only in a neighbourhood of the boundary of $\Omega$. For
$\mu\leq 1/4$ the linear problem (L) shows a remarkable structural
property for sub-harmonics, which we call Phragmen-Lindel\"of
Alternative: a given local sub-harmonic
\begin{itemize}
\item[(i)] either dominates every local super-harmonic multiplied by a suitable positive constant
\item[(ii)] or is dominated by a multiple of any local super-harmonic.
\end{itemize}
The first type of sub-harmonic is called {\em large}, the second type is
called {\em small}.

\medskip

The key to our study is the observation that solutions and
sub-solutions of \eqref{*} are sub-harmonics of (L).
We can therefore classify them according to their behavior in a neighborhood of the boundary.

A local sub-solution of \eqref{*} will be called an {\em $L$-subsolution} if it
is a large sub-harmonic and an {\em $S$-subsolution} if it is a small
sub-harmonic. In the familiar case $s=\mu=0$ large local sub-solutions are
those with finite or infinite positive boundary values and small local sub-solutions
attain zero boundary values. Note that in this paper the use of the word
``large'' for a sub-solution does not imply that this sub-solution has
``infinite boundary values''.

\smallskip

When both (L) and (N) are combined into problem (\ref{*}), interesting
thres\-hold-pheno\-mena with respect to existence or non-existence of local
sub-solutions occur. Our first main result, given in Theorem \ref{t-Main}, can be
summarized as follows: if $p>1$ and $\mu\leq 1/4$ then
\begin{center}
\fbox{
local $L$-subsolutions of \eqref{*} exist if and only if
$\displaystyle\frac{s-2}{p-1}<\frac{1}{2}-\sqrt{\frac{1}{4}-\mu}$.}
\end{center}
The proof of the main result goes as follows:
\begin{itemize}
\item[(i)] any local sub-solution $u$ of \eqref{*} satisfies the bound
  $u(x)\leq \const\delta(x)^\frac{s-2}{p-1}$, which is known as the
  Keller-Osserman upper bound
\item[(ii)] if $\mu<1/4$ then any local large sub-harmonic $u$ of (L) satisfies
$ \limsup_{x\to\partial\Omega} \frac{u(x)}{\delta(x)^{\beta_-}}>0$ where $\beta_-=\frac{1}{2}-\sqrt{\frac{1}{4}-\mu}$.
\end{itemize}
Both (i) and (ii) are compatible if
$\frac{s-2}{p-1}<\beta_-$ and incompatible
if $\frac{s-2}{p-1}>\beta_-$. The equality case belongs to the
non-existence regime, but this requires a much more refined
analysis. Likewise, the case $\mu=1/4$ is more subtle and needs
extra care.

\medskip

Our second main result, which is also given in Theorem \ref{t-Main},
shows that in the existence case, one can in fact prove the existence of two
different $L$-solutions:
\begin{itemize}
\item[(i)] an $ML$-solution to \eqref{*}, which is large but still
  dominated by at least one super-harmonic
\item[(ii)] an $XXL$-solution, which dominates every super-harmonic
  and moreover grows as fast as the Keller-Osserman upper-bound $\delta(x)^\frac{s-2}{p-1}$.
  \end{itemize}
As a consequence of the two main results we note that \eqref{*} has local sub-solutions blowing up near the boundary if and only if $s<2$ and $\mu^\ast<\mu\leq \frac{1}{4}$. Here
$\mu^\ast=\frac{1}{4}-\Big(\frac{p-2s+3}{p-1}\Big)^2$ is a negative value because $s<2$.
It is an open problem to determine the precise asymptotic behavior of an $XXL$-solution.
We conjecture that the $XXL$-solution $U(x)$ is unique and that its correct asymptotic
behaviour is given by $\lim_{x\to\partial\Omega}U(x)/\delta(x)^\frac{s-2}{p-1}=\const$
\medskip

The paper is organized as follows. In Section \ref{linear} we analyse
the linear problem (L). We explain the role played by the
Hardy-constant and prove the Phragmen-Lindel\"of
Alternative. Moreover, we construct explicit sub- and super-harmonics
and give estimates for the boundary-behaviour of large and small
sub-harmonics. In Section \ref{nonlinear} we prove a comparison
principle, which plays an important role in our analysis, and we prove
the Keller-Osserman upper bound. Section \ref{ex_nonex} contains the
proof of the main result. In Section \ref{small_conclud} we give
some additional results about small sub-solutions of \eqref{*} and in the
final Section \ref{open} we pose some open problems.

\section{Linear problem} \label{linear}
\subsection{Definitions}
For $\rho>0$ and $\eps\in(0,\rho)$ we use the notation
$$
\begin{array}{ll}
\Omega_\rho:=\{x\in\Omega:\delta(x)<\rho\}, & \Omega_{\eps,\rho}:=\{x\in\Omega:\eps<\delta(x)<\rho\} \vspace{\jot}\\
D_\rho:=\{x\in \Omega: \delta(x)>\rho\},
& \Gamma_\rho:=\{x\in\Omega:\delta(x)=\rho\}.
\end{array}
$$
\begin{figure}[ht]
\begin{center}
\scalebox{0.35}{\includegraphics[angle=270,origin=br]{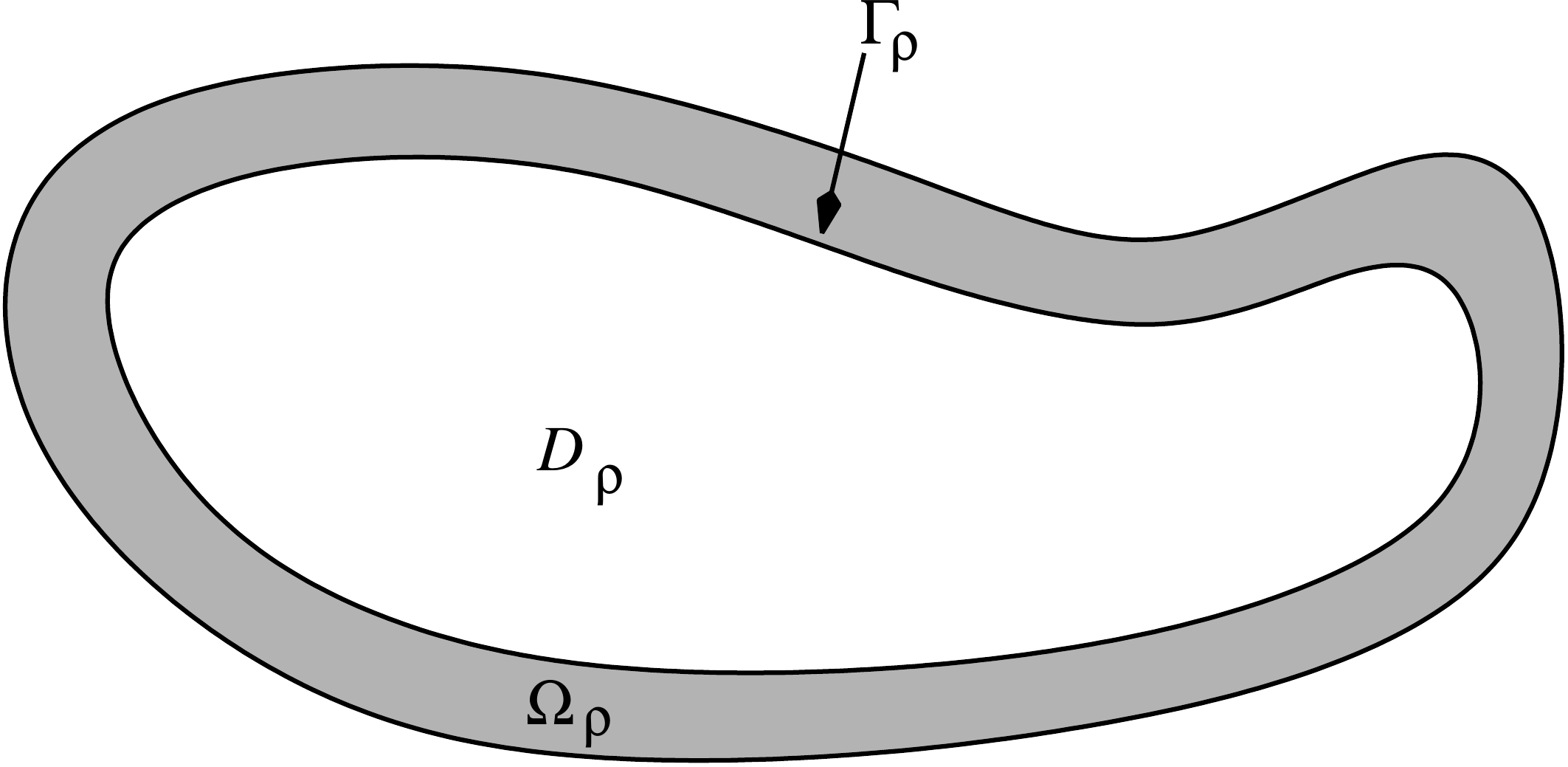}}$\quad$
\scalebox{0.35}{\includegraphics[angle=270,origin=br]{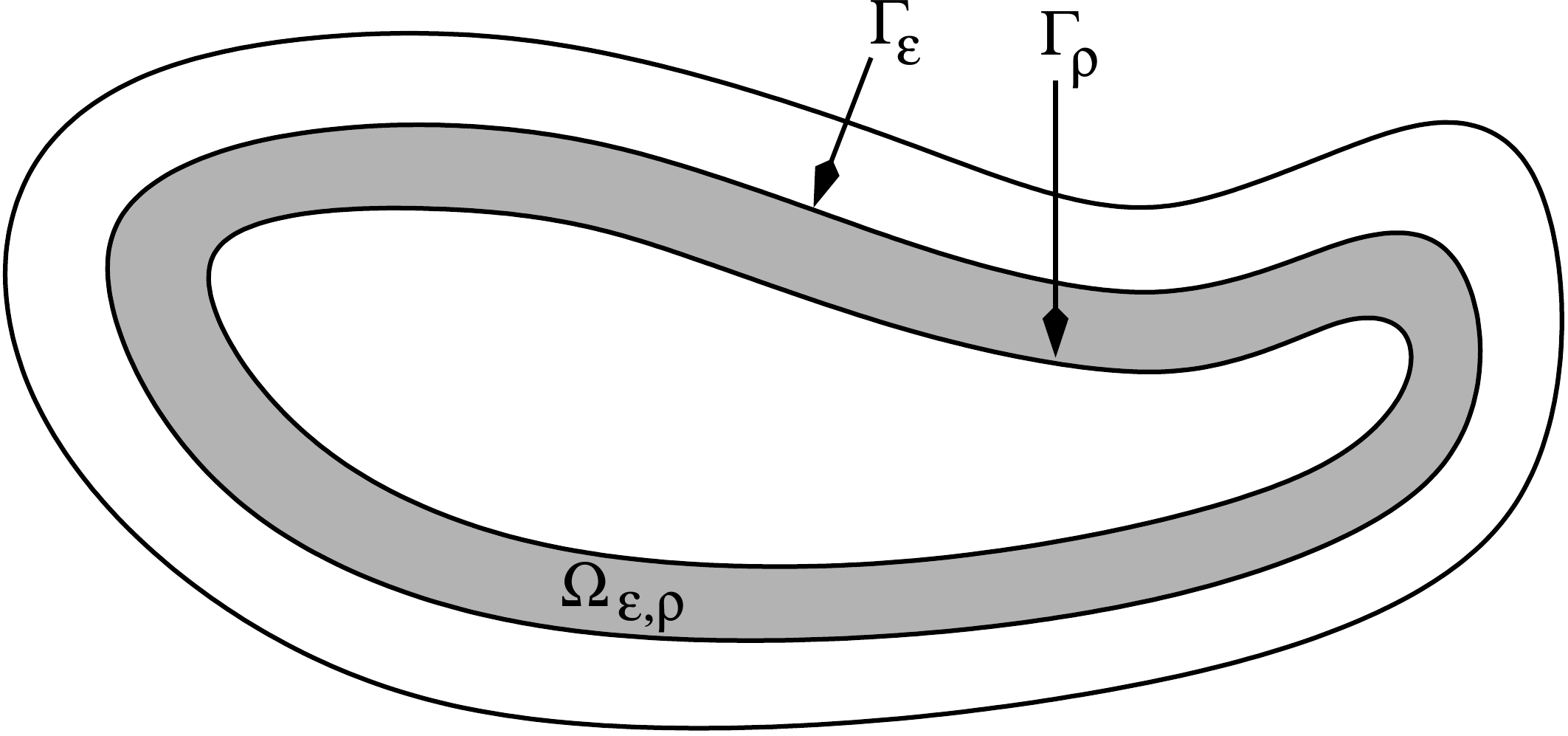}}
\end{center}
\caption{Illustration of $\Omega_\rho, D_\rho, \Gamma_\rho$ and $\Omega_{\epsilon,\rho}$}
\label{illu}
\end{figure}

In this section we present several auxiliary facts concerning the linear problem (L). For simplicity set
$$\L_\mu:=-\Delta -\frac{\mu}{\delta^2}.$$
Then (L) can be written in the form
\begin{equation}\label{e:0}
\L_\mu h=0\quad\mbox{in }\:\Omega.
\end{equation}
For convenience we call its solutions {\em harmonics}.

\begin{definition}
Let $G\subset\Omega$ and let $H^1_c(G)$ denote the space functions from $H^1(G)$ with compact
support. A sub-harmonic in $G$
is a function $\hsub\in H^1_{loc}(G)\cap C(G)$ such that
\begin{equation*}
\int_{G}\nabla \hsub\cdot\nabla\varphi\,dx-
\int_{G}\frac{\mu}{\delta^2}\hsub\varphi\,dx\leq \:0,
\qquad\forall\:0\le\varphi\in H^1_c(G).
\end{equation*}
We say that $\hsub$ is a local
sub-harmonic if there exists a parallel set $\Omega_\rho$, $\rho>0$ such that $\hsub\in
H^1_{loc}(\Omega_\rho)\cap C(\Omega_\rho)$ is a sub-harmonic in
$\Omega_\rho$. Similarly, (local) super-harmonics $\hsup$ are defined
with ``$\geq$'' in the above inequality.
\end{definition}

\begin{remark}\label{Harnack}
By the classical maximum principle for the Laplacian,
any nontrivial super-har\-monic $\hsup \gneq 0$ in $G$ is strictly positive in $G$,.
Recall also that if $\hsub$ is a sub-harmonic in $G$
then $\hsub_+$ is also a sub-harmonic in $G$, cf. \cite[Lemma 2.10]{Agmon}
\end{remark}

\subsection{The role of the Hardy constant.}

The principal result of this section is given next.
\begin{theorem} \label{t-nonexist}
Equation \eqref{e:0} admits a local positive super-harmonic
if and only if $\mu\le 1/4$. In particular no nontrivial harmonics exist if $\mu>\frac{1}{4}$.
\end{theorem}
\noindent Its proof is accomplished via the following two lemmas which are
 intimately related to Hardy's inequality. Recall that the classical
 Hardy inequality reads as follows. There exists a constant
 $C_H(\Omega)>0$ such that
\begin{equation}\label{e-Hardy}
\int_\Omega|\nabla u|^2\,dx\ge C_{H}(\Omega)\int_\Omega\frac{u^2}{\delta^2}\,dx,
\qquad\forall\:u\in H^1_0(\Omega).
\end{equation}
The optimal constant will be denoted by $C_H(\Omega)$.
For a bounded  Lipschitz domain it is known that $C_H(\Omega)\in (0,1/4]$.
If $\Omega$ is convex then $C_H(\Omega)=1/4$. In general, $C_H(\Omega)$ varies with the domain
and could be arbitrary small (see, e.g. \cite[Theorem I and Section 4]{MMP}) for a discussion and examples,
see also \cite{Davies}).

The relation between Hardy inequalities and the existence of local positive super-har\-monics
in a neighborhood of the boundary is explained by the following classical result
(cf. \cite[Theorem 3.3]{Agmon}).
\begin{lemma}\label{t-GST}
The following three statements are equivalent:
\begin{itemize}
\item[(i)] Equation \eqref{e:0} admits a positive super-harmonic in $\Omega_\rho$.
\item[(ii)] If $\hsub$ and $\hsup$ are sub- and super-harmonics of
  \eqref{e:0} in a subdomain $G$ with $\overline{G}\subset\Omega_\rho$
and if $\hsub\leq \hsup$ on $\partial G$ then $\hsub \leq \hsup$ a.e. in $G$.
\item[(iii)] The following inequality holds:
\begin{equation}\label{e-Hardy-loc}
\int_{\Omega_\rho}|\nabla u|^2\,dx\ge \mu\int_{\Omega_\rho}\frac{u^2}{\delta_\Omega^2}\,dx,
\qquad\forall\:u\in H^1_0(\Omega_\rho).
\end{equation}
\end{itemize}
\end{lemma}
Note that  the above inequality
\eqref{e-Hardy-loc} is not a particular case of \eqref{e-Hardy}
because $\dist(x,\partial\Omega)\not = \dist(x,\partial \Omega_\rho)$.
Denote the optimal constant in \eqref{e-Hardy-loc} by
$$C_H^{loc}(\Omega_\rho):=\inf_{H^1_0(\Omega_\rho)}
\frac{\int_{\Omega_\rho}|\nabla u|^2\,dx}{\int_{\Omega_\rho}\frac{u^2}{\delta^2}\,dx}.$$
The following result can be extracted from the arguments in \cite[p.3246]{MMP}.
\begin{lemma}\label{t-Hardy-loc}
{\sc (Local Hardy Inequality)}
There exists $\bar\rho>0$ such that for every $\rho\in(0,\bar\rho)$ one has $C_H^{loc}(\Omega_\rho)=1/4$.
\end{lemma}
\begin{proof}
It was already observed in \cite{BrMa97}  that $C_H^{loc}(\Omega_\rho)\ge 1/4$.
(It follows also simply from the fact that the
equation $\L_{1/4}h=0$ admits positive super-solutions in $\Omega_{\bar\rho}$ for some $\bar\rho>0$,
see Lemma \ref{l-A} below.)
On the other hand, the proof of Theorem 5 in \cite{MMP} implies that $C_H^{loc}(\Omega_\rho)\le 1/4$
for all $\rho>0$.
\end{proof}

Observe that in contrast to the "global" Hardy constant $C_H(\Omega)$ from \eqref{e-Hardy},
the value of $C_H^{loc}(\Omega_\rho)$ does not depend on the shape of domain $\Omega$
if $\rho$ is sufficiently small.

\subsection{Phragmen--Lindel\"of alternative.}

We establish a version of the Phrag\-men--Lin\-del\"of type comparison
principle for sub-harmonics, which shows that sub-harmonics are in a certain sense "separated"
by the the cone of positive super-harmonics. See \cite[pp. 93-106]{Protter} for a classical reference to
the Phragmen--Lindel\"of principle.


\begin{theorem}\label{t-PL}
{\sc (Phragmen--Lindel\"of Alternative)}
Let $\mu\le 1/4$. Let $\hsub$ be a local sub-harmonic. Then the following alternative holds:
\begin{enumerate}
\item[{$(i)$}]
either for every local super-harmonic $\hsup>0$
\begin{equation}\label{i}
\limsup_{x\to\partial\Omega}\hsub/\hsup>0,
\end{equation}
\item[{$(ii)$}]
or for every local super-harmonic $\hsup>0$
\begin{equation}\label{ii}
\limsup_{x\to \partial\Omega} \hsub/\hsup<\infty.
\end{equation}
\end{enumerate}
\end{theorem}

\begin{proof}
Assume $(i)$ does not hold, that is there exists a super-harmonic $\hsup_\ast>0$ that
\begin{equation}\label{e-A-small}
\lim_{x\to\partial\Omega} \hsub/\hsup_\ast=0.
\end{equation}
Let $\hsup>0$ be an arbitrary super-harmonic in $\Omega_\rho$.
By Remark \ref{Harnack}, there exists a constant $c>0$ such that $\hsup \ge c\hsub$ on $\Gamma_{\rho/2}$.
For $\tau>0$, define a comparison function
$$
v_\tau:=c\hsub-\tau \hsup_\ast.$$
Then \eqref{e-A-small} implies that for every $\tau>0$ there exists $\eps=\eps(\tau)\in(0,\rho)$ such that
$v_\tau\le 0$ on $\Omega_{\eps}$.
Applying the classical comparison principle in $\Omega_{\eps/2,\rho/2}$, we conclude that
$\hsup\ge v_\tau$ in $\Omega_{\eps/2,\rho/2}$ and hence, in $\Omega_{\rho/2}$.
So by considering arbitrary small $\tau>0$, we conclude that
for every super-harmonic $\hsup>0$ in $\Omega_\rho$ there exist $c>0$
such that $\hsup\geq c\hsub$ holds in $\Omega_\rho$. This implies \eqref{ii}.
\end{proof}

Theorem \ref{t-PL} suggests the following classification of sub-harmonics.
\begin{definition}\label{d-SL}
Let $\mu\le 1/4$ and let $\hsub$ be a local sub-harmonic in $\Omega_\rho$.
We say that $\hsub$ is {\em large} if it satisfies the first alternative (i).
Otherwise, we say that $\hsub$ is a {\em small}.
\end{definition}

The classification of harmonics into small and large harmonics is
included in the above definition.
In the sequel we shall use the notation $\hsub$ for small and $\Hsub$ for large sub-harmonics.

\subsection{Construction of local sub- and super-harmonics.}
It is well known (cf. \cite[Lemma 14.15]{Gilbarg}) that if $\Omega$
is of class $C^k$, $k\ge 2$, then there exists $\bar\rho>0$ such that
the distance function $\delta$ is in $C^k(\Omega_{\bar\rho})$ and the
set $\Gamma_\eps$ is of class $C^k$ for all $\eps\in(0,\bar\rho)$.
For every $x\in\Gamma_\eps$ there exists a unique point $\sigma(x)\in\partial\Omega$
such that $|x-\sigma(x)|=\delta(x)$. Moreover,
\begin{equation}\label{e-nabla}
|\nabla\delta(x)|=1+o(\delta(x))\quad\text{as }\:\delta(x)\to 0,
\end{equation}
while
\begin{equation}\label{e-delta}
\Delta\delta(x)=-(N-1)\H_0(\sigma(x))+o(\delta(x))\quad\text{as }\:\delta(x)\to 0,
\end{equation}
where $\H_0(\sigma(x))$ denotes the mean curvature of $\partial\Omega$ at the point $\sigma(x)$.
Note that the mean curvature of $\partial\Omega$ is bounded, since $\Omega$ is bounded and smooth.
\medskip

In what follows,
$\beta_-\le \beta_+$ denote the real roots of the scalar equation
$\beta(1-\beta)=\mu$, i.e.
\beqn \label{roots}
 \beta_{\pm} = \frac{1}{2}\pm
\sqrt{\frac{1}{4}-\mu}
\tx{ provided } \mu\le 1/4.
\eeqn
Clearly, if $\mu=1/4$ then $\beta_-=\beta_+=1/2$.
\begin{lemma}\label{l-A}
$(i)$
Let $\mu<1/4$.
The function $\delta^\beta$ is a local super-harmonic of $\L_\mu$ if $\beta\in(\beta_-,\beta_+)$
and a local sub-harmonic of $\L_\mu$ if $\beta\not\in[\beta_-,\beta_+]$.
Moreover, if $\epsilon\in(0,\epsilon_\ast)$, where  $\epsilon_\ast=\min\{1,\sqrt{1-4\mu}\}$
then
$$
\hsup:=\delta^{\beta_+}(1-\delta^{\epsilon}),\qquad
\Hsup:=\delta^{\beta_-}(1+\delta^{\epsilon})
$$
are positive local super-harmonics of $\L_\mu$,
while
$$
\hsub:=\delta^{\beta_+}(1+\delta^{\epsilon}),\qquad
\Hsub:=\delta^{\beta_-}(1-\delta^{\epsilon})
$$
are positive local sub-harmonics of $\L_\mu$.

$(ii)$
Let $\mu=1/4$.
The function $\delta^{1/2}\left(\log\frac{1}{\delta}\right)^\beta$ is
a local super-harmonic of $\L_{1/4}$ if $\beta\in(0,1)$
and a local sub-harmonic of $\L_{1/4}$ if $\beta\not\in[0,1]$.
Moreover, if $\epsilon\in(0,1)$ then
$$
\hsup=\delta^{1/2}\left(1-\Big(\log\frac{1}{\delta}\Big)^{-\epsilon}\right),\qquad
\Hsup:=\delta^{1/2}\log\frac{1}{\delta}\left(1+\Big(\log\frac{1}{\delta}\Big)^{-\epsilon}\right)
$$
are positive local super-harmonics of $\L_{1/4}$, while
$$
\hsub:=\delta^{1/2}\left(1+\Big(\log\frac{1}{\delta}\Big)^{-\epsilon}\right),\qquad
\Hsub:=\delta^{1/2}\log\frac{1}{\delta}\left(1-\Big(\log\frac{1}{\delta}\Big)^{-\epsilon}\right)
$$
are positive local sub-harmonics of $\L_{1/4}$.
\end{lemma}

\begin{proof}
$(i)$
Note that
\begin{eqnarray*}
\nabla\delta^\beta&=&
\beta\delta^{\beta-1}\nabla\delta,\\
-\Delta\delta^\beta&=&
\beta(1-\beta)\delta^{\beta-2}|\nabla\delta|^2-\beta\delta^{\beta-1}\Delta\delta.
\end{eqnarray*}
Thus a direct computation together with \eqref{e-nabla}, \eqref{e-delta} imply the result
(cf. \cite[Lemma 7]{MMP}).

$(ii)$
Observe that
\begin{eqnarray*}
\nabla\left(\delta^{\gamma}\log^\beta\frac{1}{\delta}\right)&=&
\left(\gamma\log^\beta\frac{1}{\delta}-\beta\log^{\beta-1}\frac{1}{\delta}\right)
\delta^{\gamma-1}\nabla\delta,\\
-\Delta\left(\delta^{\gamma}\log^\beta\frac{1}{\delta}\right)&=&
\left(\gamma(1\!-\!\gamma)\log^\beta\frac{1}{\delta}
+\beta(2\gamma\!-\!1)\log^{\beta-1}\frac{1}{\delta}
+\beta(1\!-\!\beta)\log^{\beta-2}\frac{1}{\delta}\right)
\delta^{\gamma-2}|\nabla\delta|^2\\
&-&
\left(\gamma\log^\beta\frac{1}{\delta}-\beta\log^{\beta-1}\frac{1}{\delta}\right)
\delta^{\gamma-1}\Delta\delta.
\end{eqnarray*}
Thus a direct computation together with \eqref{e-nabla}, \eqref{e-delta} imply the result.
\end{proof}
The following theorem,
which is an immediate corollary of Theorem \ref{t-PL} and Lemma \ref{l-A},
summarises our results concerning the asymptotic behaviour of sub-harmonics at the boundary.

\begin{theorem}\label{table}
Let $\hsub$ be a small local sub-harmonic and $\Hsub$ be a large local sub-harmonic of $\L_\mu$.
\begin{itemize}
\item[(i)] If $\mu<1/4$ then
$$\limsup_{x\to \partial\Omega}\frac{\hsub}{\delta^{\beta_+}}<\infty,\qquad
\limsup_{x\to \partial\Omega}\frac{\Hsub}{\delta^{\beta_-}}>0.$$
\item[(ii)] If $\mu=1/4$ then
$$\limsup_{x\to \partial\Omega}\frac{\hsub}{\delta^{1/2}}<\infty,\qquad
\limsup_{x\to
  \partial\Omega}\frac{\Hsub}{\delta^{1/2}\log\frac{1}{\delta}}>0.$$
\end{itemize}
The above leading order terms are sharp.
\end{theorem}


Lemma \ref{l-A} and Theorem \ref{table} have the following implications.
\begin{corollary} Let $\mu\leq 1/4$.
\begin{itemize}
\item[(i)] The small local sub-harmonics vanish on the entire boundary of $\Omega$.
\item[(ii)] If $\mu<0$ then the large sub-harmonics are unbounded at some points of
  $\partial \Omega$.
\item[(iii)] If $0<\mu\le 1/4$ then there exist large sub-harmonics vanishing on $\partial
\Omega$.
\end{itemize}
\end{corollary}

\begin{remark}
\begin{enumerate}
\item Observe that when $\mu=0$ then $\Hsub=const.$ is a large sub-harmonic in a neighbourhood of the boundary.
\item For $\mu<1/4$ large local sub-harmonics fail to belong to the
  subspace of functions in $H^1(\Omega_\rho)$ which vanish on
  $\partial\Omega$. Indeed, for $\mu\leq 0$ large local sub-harmonics do not
  converge to zero near $\partial\Omega$. And for $0<\mu<1/4$, even if
  a large local sub-harmonic vanishes on $\partial\Omega$ then its
  gradient is not square-integrable near $\partial\Omega$. To
see this, let $\Hsub$ be a large sub-harmonic of $\L_\mu$ in $\Omega_\rho$.
For $\beta\in (1/2,\beta_+)$ the function $\delta^\beta$ is a
super-harmonic in $H^1(\Omega_\rho)$ with vanishes on
$\partial\Omega$. Hence $\Hsub_\kappa:=(\Hsub-\kappa\delta^\beta)_+$ is a large sub-harmonic.
By choosing a sufficiently large $\kappa>0$ we can ensure that $\Hsub_\kappa$ vanishes
on $\Gamma_\rho$. Assume for contradiction that $\Hsub_\kappa\in H^1_0(\Omega_\rho)$. Then
$$\int_{\Omega_\rho}|\nabla \Hsub_\kappa|^2\,dx-\mu\int_{\Omega_\rho}\frac{\Hsub_\kappa^2}{\delta_\Omega^2}\,dx\le 0$$
and by the local Hardy inequality we obtain $\Hsub_\kappa\equiv 0$, i.e.,
$\Hsub\leq \kappa\delta^\beta$. This contradicts Theorem \ref{table}(i).

\item If $\mu>1/4$,
then $\L_\mu$ has no positive local super-harmonics (cf. Theorem~\ref{t-nonexist}). However,
$$
\delta^{\gamma}\log^\beta\frac{1}{\delta}
$$
is a local sub-harmonic for arbitrary $\gamma,\beta\in\R$.
This suggests that in the case $\mu>1/4$ local
sub-harmonics can not be naturally classified according to their asymptotic behaviour.
\end{enumerate}
\end{remark}

Another direct consequence of Theorem \ref{t-PL} and Lemma \ref{l-A}
is a two--sided bound on the asymptotic behaviour of positive super-harmonics at the boundary.

\begin{theorem}\label{table-sup}
Let $\hsup>0$ be a local super-harmonic of $\L_\mu$.
\begin{itemize}
\item[(i)] If $\mu<1/4$ then
$$\liminf_{x\to \partial\Omega}\frac{\hsup}{\delta^{\beta_+}}>0,\qquad
\liminf_{x\to \partial\Omega}\frac{\hsup}{\delta^{\beta_-}}<\infty.$$
\item[(ii)] If $\mu=1/4$ then
$$\liminf_{x\to \partial\Omega}\frac{\hsup}{\delta^{1/2}}>0,\qquad
\liminf_{x\to \partial\Omega}\frac{\hsup}{\delta^{1/2}\log\frac{1}{\delta}}<\infty.$$
\end{itemize}
The above leading order terms are sharp.
\end{theorem}




\section{Estimates for the nonlinear problem}\label{nonlinear}
\subsection{Comparison principle}
We start with the definition of sub- and super-solutions to the
nonlinear problem \eqref{*}.
\begin{definition}
A \textit{sub-solution} to \eqref{*} in a subdomain
$G\subset\Omega$ is a function $\usub \in H^1_{loc}(G)\cap C(G)$ such that
\begin{equation}
\label{weak}
\int_{G}\nabla \usub \cdot\nabla\varphi\,dx-
\int_{G}\frac{\mu}{\delta^2}\usub \varphi\,dx\,
+\int_{G}\frac{\usub^p}{\delta^s}\varphi\,dx\le 0,
\quad\forall\:0\le\varphi\in H^1_c(G).
\end{equation}
A {\em super-solution} $\usup$ is defined similarly
by replacing "$\le$" with "$\ge$". A function $u$ which is both a sub- and super-solution will be called a {\em solution}.
\end{definition}
\begin{lemma} {\sc (Comparison principle)} \label{comparison}
\begin{enumerate}
\item[$(i)$] Let $G$ be open with $G\subset\Omega$.
Let $0\le \usub, \usup \in H^1_{loc}(G)\cap C(G)$ be a pair of sub- and super-solutions to \eqref{*} in $G$
such that $$\limsup_{x\to\partial G}[ \usub(x)-\usup(x)]<0.$$
Then $\usub \leq \usup$ in $G$.
\item[$(ii)$]
Let $G$ be open with $\overline{G}\subset\Omega$.
Let $\usub, \usup\in H^1(G) \cap C(\overline{G})$ be a pair of sub- and super-solutions
to \eqref{*} in $G$ with $\usup>0$ in $G$ and $\usub\leq \usup$ on $\partial G$.
Then $\usub\leq \usup$ in $G$.
\label{hardy_cp}
\end{enumerate}
\end{lemma}
\begin{proof}
$(i)$
Subtracting one inequality from another we obtain

\begin{align*}
\int_{G}\nabla( \usub- \usup)\cdot\nabla\varphi\,dx-
\int_{G}\frac{\mu}{\delta^2}(\usub-\usup) \varphi\,dx\,
+\int_{G}\frac{W(x)}{\delta^s}(\usub-\usup)\varphi\,dx\le 0,\\
\quad\forall\varphi\in H^1_c(G),\varphi\geq 0,
\end{align*}
where
$$W(x):=\frac{\ul u^p-\ol u^p}{\ul u-\ol u}.$$
Assume that $(\ul u-\ol u)_+\not\equiv 0$. Testing against $(\ul u-\ol u)_+$ we conclude that
\begin{align}\label{ineq1}
\int_G\left(|\nabla(\ul u-\ol u)_+|^2-\frac{\mu}{\delta^2}(\ul u-\ol
  u)_+^2+\frac{W(x)}{\delta^s}(\ul u-\ol u)_+^2\right)\:dx\le 0.
\end{align}
Since $\ol u>0$ we can write
$$
(\usub-\usup)_+=\usup\phi,
$$
where $\phi\in H^1_c(G)$ due to the assumption that
$\limsup_{x\to\partial G}[ \usub(x)-\usup(x)]<0$. Note that
$\supp\phi=\overline{G}_+$, where $G_+:= \{x\in G:
\usub(x)>\usup(x)\}$. We obtain
\begin{align*}
\int_G |\nabla(\usub-\usup)_+|^2\:dx&= \int_G(\phi^2|\nabla \usup|^2+2\phi\usup \nabla \usup\cdot \nabla \phi
+\usup^2 |\nabla \phi|^2)\:dx\\
&=\int_G\phi^2|\nabla \usup|^2\:dx +\int_G\nabla \usup\cdot
\nabla(\phi^2\usup) \:dx\\
& \geq\int_G \frac{\mu}{\delta^2}\phi^2 \usup^2 - \frac{\phi^2\usup^{p+1}}{\delta^s}\:dx,
\end{align*}
where we have used that $\usup$ is a super-solution. Hence we conclude that
\begin{align}\label{ineq2}
\int_G\left(|\nabla(\ul u-\ol u)_+|^2-\frac{\mu}{\delta^2}(\ul u-\ol u)_+^2+\frac{V(x)}{\delta^s}(\ul u-\ol u)_+^2\right)\ge 0,
\end{align}
where $V(x)=\ol u^{p-1}$.
But by strict convexity we have $W(x)\gneq V(x)$ on $G_+$. Thus
\eqref{ineq2} and \eqref{ineq1} imply that $G_+$ has zero measure,
which contradicts the assumption $(\usub-\usup)_+\not\equiv 0$.

The proof of $(ii)$ is similar if instead of $\phi$ one uses $\phi_\eps$
defined by $(\usub-(\usup+\eps))_+=\usup\phi_\eps$ with $\eps>0$, so we omit it.
\end{proof}

\begin{remark}
Note that the above lemma is valid for any $\mu\in\R$.
We do not require the assumption $\mu\le1/4$ which ensures positivity
of the principal part $\L_\mu$ because for $\mu>1/4$ the nonlinearity
compensates for the loss of positivity.
\end{remark}

\subsection{Keller--Osserman type bound}
\label{kob}

By a simple computation analogous to Lemma \ref{l-A} one finds that
for $p>1$ the function
$$
\gamma \delta(x)^\frac{2-s}{1-p}
$$
has the following properties:

\begin{table}[h!]
\begin{center}
\begin{tabular}{|c|c|c|}
\hline
 & local sub-solution & local super-solution \\ \hline
$\beta_-\leq \frac{2-s}{1-p}\leq\beta_+$ & --  & $\gamma$ arbitrary \\ \hline
$\frac{2-s}{1-p}<\beta_-$ or $\beta_+<\frac{2-s}{1-p}$ & $\gamma$ small
& $\gamma$ large \\ \hline
\end{tabular}
\end{center}
\caption{Properties of $\delta(x)^\frac{2-s}{1-p}$}
\label{ko_tab}
\end{table}
In particular, this function is always a local super-solution if
$\gamma$ is sufficiently large. The next considerations show that in
order to make it a global super-solution, one needs to replace the
distance function $\delta$ by the regularized distance function $d:\Omega \to \mathbb{R}^+$
attributed to Whitney, cf. \cite{S}.
The regularized distance function is in $C^\infty(\Omega)$ regardless
of the regularity of $\partial \Omega$ and has the following
properties: there exists a positive constant $c$ such that
\beqn\label{dist}
c^{-1}\delta(x)\leq d(x)&\leq& c\delta(x),\\
\nonumber |\nabla d(x)|&\leq& c,\\
\nonumber |\Delta d(x)|&\leq& cd^{-1}(x)\quad\tx{for all }\: x\in \Omega.
\eeqn
\begin{proposition}\label{supersol}
Let $p>1$. For $\gamma$ sufficiently large, but independent of $\eps\geq 0$,
the function
$$
\usup = \gamma d^{\frac{s}{p-1}}(d-\eps)^{-\frac{2}{p-1}}
$$
is a super-solution  of \eqref{*} in $\{x\in \Omega, d(x)>\eps\}$.
\end{proposition}
\begin{proof}
A straightforward computation together with \eqref{dist} yields
\beq
|\Delta  d^{\frac{s}{p-1}}| &\leq& \alpha_1  d^{\frac{s}{p-1}-2} \;\leq \;\alpha_1 d^\frac{s}{p-1}(d-\epsilon)^{-2},\\
|\Delta (d-\eps)^{-\frac{2}{p-1}}| &\leq& \alpha_2
(d-\eps)^{-\frac{2}{p-1}-2}+\alpha_3(d-\eps)^{-\frac{2}{p-1}-1}d^{-1}\\
&\leq&(\alpha_2+\alpha_3)(d-\epsilon)^{-\frac{2p}{p-1}},\\
|(\nabla  d^{\frac{s}{p-1}}, \nabla (d-\eps)^{-\frac{2}{p-1}})|&\leq&
\alpha_4 d^{\frac{s}{p-1}-1}(d-\eps)^{-\frac{2}{p-1}-1}\;\leq\;\alpha_4d^\frac{s}{p-1}(d-\epsilon)^{-\frac{2p}{p-1}},
\eeq
where $\alpha_i$, $i=1\dots 4$ depend only on $c,p,s$. In addition
$$
|\frac{\mu}{\delta^2}\usup| \leq \alpha_5\gamma d^{\frac{s}{p-1}}(d-\eps)^{\frac{-2p}{p-1}},
$$
where again $\alpha_5$ depends on $c,\mu, p,s$. Collecting all the terms and keeping in mind
that
$$
\delta^{-s}(x)\geq c^{-|s|}d^{-s}(x)
$$
we find
\beq
\L_\mu \usup &\geq& -\gamma \alpha_6 d^{\frac{s}{p-1}}(d-\eps)^{-\frac{2p}{p-1}}
\geq -\gamma\alpha_6 c^{|s|} d^{\frac{s}{p-1}+s}(d-\eps)^{-\frac{2p}{p-1}}\delta^{-s}\\
&=&- \alpha_6c^{|s|}\frac{\usup^p}{\gamma^{p-1}}\delta^{-s} \geq- \usup^p\delta^{-s},
\eeq
for $\gamma$ sufficiently large, but independent of $\eps\geq 0$.
\end{proof}

Sub-solutions to the nonlinear equation \eqref{*} obey a universal upper bound given next.
As a tool we use the comparison principle from Lemma \ref{comparison}.
\begin{proposition}
\label{ko}
{\sc (Keller--Osserman Bound)}
Assume $p>1$. Let $u$ be an arbitrary local sub-solution to \eqref{*} in $\Omega_\rho$ for some $\rho>0$.
Then there exists $\gamma_\ast>0$ depending on $u$ such that
\begin{equation}\label{KO}
u(x)\le \gamma_\ast\delta^\frac{2-s}{1-p}(x)\quad\mbox{in }\:\Omega_\rho.
\end{equation}
If $u$ is sub-solution in all of $\Omega$, then $\gamma_\ast$ can be
chosen independently of $u$.
\end{proposition}
\begin{proof}
Let $u$ be a local sub-solution of \eqref{*} in $\Omega_\rho$. Thus
$$
u(x)\leq \usup (x)= \gamma d^{\frac{s}{p-1}}(x)(d(x)-\eps)^{-\frac{2}{p-1}}\quad\tx{in}\quad\{x\in \Omega_\rho: d(x)>\eps\},
$$
provided $0<\epsilon<\rho/c$ with $c$ as in \eqref{dist} and provided
$\gamma$ is so large that $\usup\geq u$ on $\Gamma_\rho$. Since the
above inequality holds for arbitrary positive $\eps<\rho/c$ it follows that
$$
u(x)\leq \gamma d^{\frac{s-2}{p-1}}(x)\leq \gamma c^{\frac{s-2}{p-1}}\delta^{\frac{s-2}{p-1}}(x),
$$
as required. If $u$ is a sub-solution in all of $\Omega$ then the above
construction works on the set $\{x\in\Omega : d(x)>\eps\}$, which has
only the boundary at $d(x)=\eps$ and no second boundary $\Gamma_\rho$.
\end{proof}


\section{The main results}\label{ex_nonex}
Since every solution and sub-solution of \eqref{*} is a sub-harmonic of $\L_\mu$,
we shall classify them in accordance with Definition \ref{d-SL}.

\begin{definition}\label{d-sol}
A solution of \eqref{*} is called an {\it S--solution} if it is a
small sub-harmonic and it is called an {\it L--solution}
if it is is a large sub-harmonic. Further, we introduce
different classes of $L$--solutions:
\begin{description}
\item[$(ML)$]
$U$ is an {\it $ML$-solution}\footnote{{\em Moderate solutions}, as introduced in \cite{Dynkin}}
if there exists a super-harmonic $H$ such that
$$\limsup_{x\to\partial\Omega}\frac{U}{H}<+\infty;$$
\item[$(XL)$]
$U$ is an {\it $XL$--solution} of \eqref{*} if for every super-harmonic $H$ one has
$$\liminf_{x\to\partial\Omega}\frac{U}{H}=+\infty;$$
\item[$(XXL)$]
$U$ is an {\it $XXL$--solution} of \eqref{*} if one has
$$\liminf_{x\to\partial\Omega}\frac{U}{\delta^{\frac{s-2}{p-1}}}>0.$$
\end{description}
The corresponding classes of sub-solutions and local (sub)\,solutions are defined accordingly.
\end{definition}

\begin{remark}
Note that division of $L$-solutions into $ML$, $XL$, $XXL$ solutions is not exhaustive.
For example, the solution of the problem
$$\left\{
\begin{array}{l}
-\Delta u + u^p = 0
\quad\text{in $\Omega$},\\
\text{$u = 0$ on $\Gamma_0$},\quad \text{$u=1$ on $\Gamma_1$},
\quad\text{$u=+\infty$ on $\partial\Omega\setminus(\Gamma_0\cup\Gamma_1$)},\\
\end{array}
\right.
$$
where $\Gamma_0,\Gamma_1\subset\partial\Omega$ are smooth submanifolds of $\partial\Omega$,
is an $L$-solution which does not belong to the classes $ML$, $XL$, $XXL$.
\end{remark}

\begin{figure}[ht]
\begin{center}
\scalebox{1.1}{\includegraphics[angle=0,origin=br]{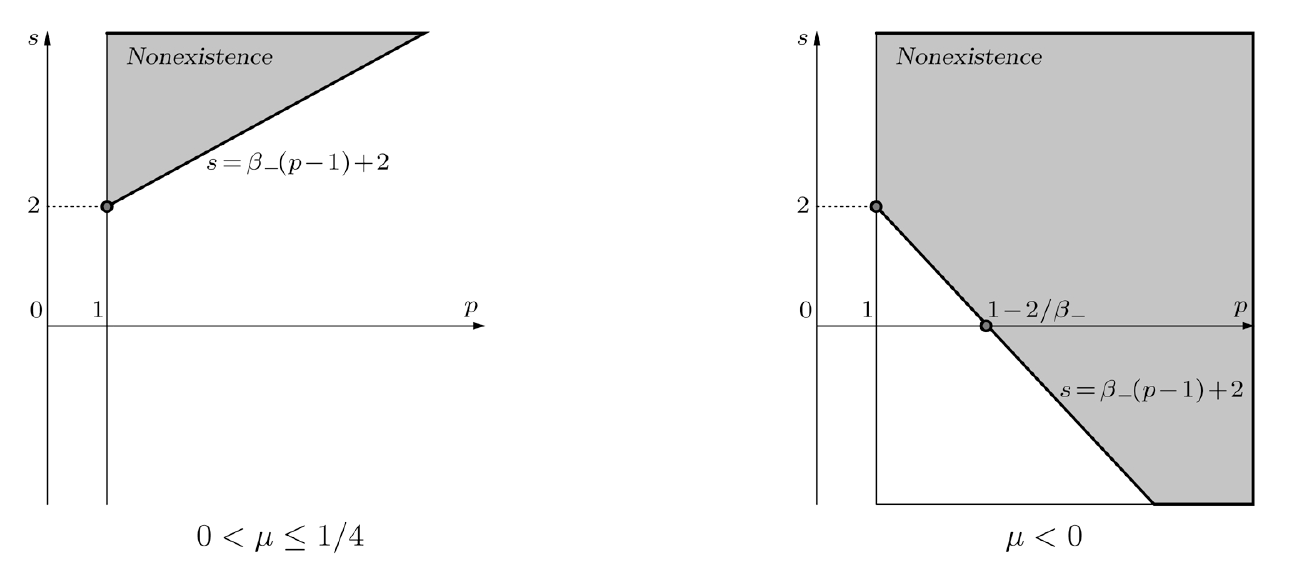}}
\end{center}
\caption{Nonexistence zones of equation \eqref{*} for typical values of $\mu$.}
\label{pic1}
\end{figure}

Our main result in the paper reads as follows.

\begin{theorem}\label{t-Main}
Let $\mu\le 1/4$, $\beta_-$ be as in \eqref{roots} and $p>1$.
\begin{itemize}
\item[i)]
If $s \geq \beta_-(p-1)+2$ then \eqref{*} admits no local $L$-subsolutions.
\item[ii)]
If $s < \beta_-(p-1)+2$ then \eqref{*} admits $ML$- and $XXL$-solutions in $\Omega$.
\end{itemize}
\end{theorem}

The above result can be seen as a critical threshold phenomenon in two
different ways by either taking $p$ or $\mu$ as a parameter.
\begin{itemize}
\item[(a)] {\em Critical value of $p$:} Let $p^\ast=1-\frac{2-s}{\beta_-}$
  with the convention $p^\ast=+\infty$ if $\beta_-=0$ and $s<2$,
  $p^\ast=-\infty$ if $\beta_-=0$ and $s\geq 2$.
\noindent
\begin{center}
\begin{tabular}{|c|c|l|}
\hline
existence & nonexistence & \\ \hline
$p\geq p^\ast$ & $1<p<p^\ast$ & if $\mu<0$\\ \hline
$1<p\leq p^\ast$ & $p>\max\{1,p^\ast\}$ & if $0<\mu\leq 1/4$\\ \hline
\end{tabular}
\end{center}
\item[(b)] {\em Critical value of $\mu$:} Let
$\mu^\ast=\frac{1}{4}-\Big(\frac{p-2s+3}{p-1}\Big)^2$.
\noindent
\begin{center}
\begin{tabular}{|c|c|l|}
\hline
existence & nonexistence & \\ \hline
---& $\mu\leq1/4$ & if $s\geq (p+3)/2$ \\ \hline
$\mu^\ast<\mu\leq 1/4$ & $\mu\leq\mu^\ast$ & if $s<(p+3)/2$\\ \hline
\end{tabular}
\end{center}
\end{itemize}

\begin{remark}
If $s\ge 2$, $\mu=0$, $p=\frac{N+2}{N-2}$ and if $\Omega$ is the unit ball in $\R^N$,
$N\ge 3$, Ratto et al. \cite{Veron} proved that no global positive solution exists.
Since such solutions are $L$-solutions our result (i) extends the non-existence
result in \cite{Veron}.
\end{remark}



In the remaining part of this section we prove Theorem \ref{t-Main}.
First we present the nonexistence part of the proof and
after that, we consider the existence.

\subsection{Proof of Theorem \ref{t-Main} }
\subsubsection { Nonexistence}

Observe that in the {\em supercritical case} $s>\beta_-(p-1)+2$
the Keller--Osserman bound \eqref{KO} is incompatible with the lower bound
on large sub-harmonics in Theorem \ref{table}. As every $L$--subsolution
to \eqref{e:0} is a large sub-harmonic of $\L_\mu$, this immediately
implies the nonexistence of local $L$--subsolutions to \eqref{e:0}.


In the {\em critical case} $s=\beta_-(p-1)+2$ the Keller--Osserman bound
is comparable with the lower bound on large sub-harmonics,
so different arguments must be used to prove the nonexistence.

Below we present a proof which covers both subcritical and critical
cases. It consists of three parts:
\begin{itemize}
\item[(a)] First we show that for every local $L$-subsolution $\usub$
  there exists a local $L$-subsolution $\usub_\ast$, which vanishes on
  $\Gamma_\rho$ and satisfies $\limsup_{x\to\partial\Omega}\frac{\usub_{\,\ast}}{\usub}=1$.
\item[(b)] Then we construct a family of super-solutions $u_\epsilon$
  in $\Omega_{\rho,\epsilon_\rho}$, converging to zero as $\epsilon\to
  0$ and tending to $+\infty$ on the inner boundary and to zero on the
  outer boundary, cf. Figure \ref{illu3}.
\begin{figure}[ht]
\begin{center}
\scalebox{0.45}{\includegraphics[angle=270,origin=br]{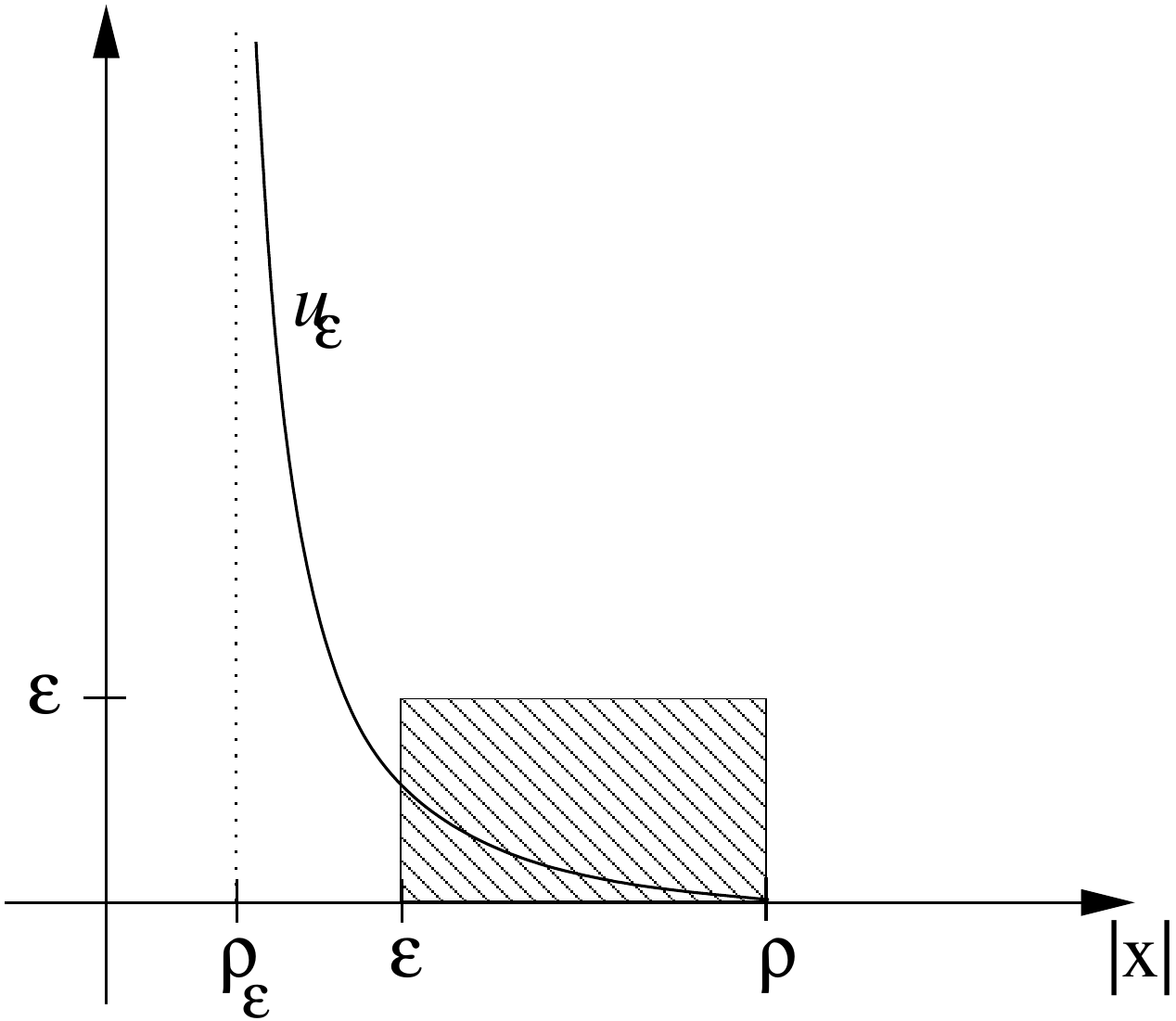}}
\end{center}
\caption{Graph of $u_\epsilon$}
\label{illu3}
\end{figure}
\item[(c)] From this "improved upper bound" it follows from by
  Comparison Principle of Lemma \ref{comparison} that zero is the only
  local $L$-subsolution.
\end{itemize}

\begin{lemma}\label{r-vanish}
Let $\mu\le 1/4$ and $p>1$.
Let $\usub\ge 0$ be a local $L$--subsolution to \eqref{*}.
Then \eqref{*} admits a local $L$--subsolution $\usub_{\,\ast}$ such that for some $\rho>0$
\begin{equation}
\usub_{\,\ast}=0\quad\text{on }\:\Gamma_\rho\quad\text{and}\quad
\limsup_{x\to\partial\Omega}\frac{\usub_{\,\ast}}{\usub}=1.
\end{equation}
\end{lemma}

\begin{proof}
Let $\usub$ be a local $L$--subsolution to \eqref{*}.
For $\kappa>0$, set
$$h_\ast:=\kappa\delta^{1/2}\log(1/\delta)^{1/2}.$$
According to Lemma \ref{l-A}, $h_\ast$ is a local super-harmonic for
$\L_{1/4}$ and hence for $\L_\mu$.
Then $\usub-h_\ast$ is also a local sub-harmonic of $\L_\mu$ and
$$\limsup_{x\to\partial\Omega}\frac{\usub-h_\ast}{h_\ast}=+\infty,$$
by Theorem \ref{table}.
In particular, this means that $\usub-h_\ast$ is a large sub-harmonic of $\L_\mu$,
according to Definition \ref{d-SL}.
Moreover, for a sufficiently small $\rho>0$ we ensure $\usub-h_\ast\le 0$ on $\Gamma_\rho$
by choosing $\kappa>0$ sufficiently big.
Besides being a large local sub-harmonic the function $(\usub-h_\ast)_+$ satisfies
\begin{equation}\label{a-sub}
\L_\mu (\usub-h_\ast)_+ \leq \L_\mu\usub \leq
-\frac{C}{\delta^s}\usub^p\leq -\frac{C}{\delta^s}(\usub-h_\ast)^p
\quad\mbox{in } \{x\in \Omega_\rho: \usub-h_\ast>0\}.
\end{equation}
Thus, setting $\usub_\ast=(\usub-h_\ast)_+$ we obtain a local
$L$--subsolution in $\Omega_\rho$ with the required properties. Note
that we have used the fact that the maximum of two sub-solutions is
again a sub-solution.
\end{proof}

\begin{remark}
Note that the function $\usub_{\,\ast}$ extended by zero to $D_\rho$ is a sub-solution to \eqref{*}
in the entire domain $\Omega$.
\end{remark}

Now we establish an "improved upper bound" on local sub-solutions of \eqref{*}
vanishing on $\Gamma_\rho$ for some $\rho>0$, which immediately implies
Theorem \ref{t-Main} $(i)$ via the Comparison Principle (Lemma \ref{comparison})
and Lemma \ref{r-vanish}.

\begin{lemma}\label{KO-impr}
Let $\mu\le 1/4$, $p>1$ and $s\ge\beta_-(p-1)+2$.
Then there exists $\rho>0$ such that for every $\eps\in(0,\rho)$
there exists $\rho_\eps\in(0,\eps)$ and
a positive super-solution $z_\eps$ of the nonlinear equation \eqref{*}
in the ring-shaped domain $\Omega_{\rho_\eps,\rho}$
such that
\begin{equation}\label{barrier}
\text{$z_\eps=0$ on $\Gamma_\rho$,\qquad
$\liminf_{x\to\Gamma_{\rho_\eps}}z_\eps=+\infty$\quad and\quad
$\sup_{\Omega_{\eps,\rho}}z_\eps\le\eps$.}
\end{equation}
\end{lemma}

\begin{remark}
In fact, Lemma \ref{KO-impr} implies more then mere nonexistence.
Consider a family of "large solution" problems
\begin{equation}\label{*-rho}
-\Delta u-\frac{\mu}{\delta^2}u+\frac{u^p}{\delta^s}=0\quad\mbox{in } D_\eps,\qquad
u=+\infty \quad\mbox{on } D_\eps.
\end{equation}
For each $\rho>0$ such a problem is well-posed and admits a unique "large solution" $u_\eps$,
cf. \cite{MV2}. Moreover, the family $u_\eps$ is monotone nonincreasing on compact subsets of $\Omega$
as $\eps\to 0$. Thus for all sufficiently small $\eps>0$ one can find a super-harmonic $h$ in $\Omega_\rho$,
where $\rho>0$ is taken from Lemma~\ref{KO-impr}, so that $u_\eps-h\le 0$ on $\Gamma_\rho$.
If $s\ge\beta_-(p-1)+2$ then Lemma~\ref{KO-impr} implies that
$(u_\eps-h)_+$ converges to zero as $\eps\to 0$, uniformly on every compact subset of $\Omega$.
Thus, in the nonexistence regime $s\ge\beta_-(p-1)+2$,
an attempt to approximate solutions of \eqref{*} by exhausting the domain
$\Omega$ will lead to an $S$--solution (possibly trivial) in the limit.
\end{remark}


\paragraph{{\it Proof of Lemma \ref{KO-impr}.}}
We are going to construct the super-solutions $u_\eps$ satisfying \eqref{barrier}
using the solutions of an ODE initial value problem.
Related arguments were previously used in \cite{KLMS}.

Let as before $\sigma(x)$ be the projection of the point $x\in\Omega_\rho$ on $\partial\Omega$
and $\delta(x)$ be the distance of $x$ to the boundary.
Fix $\bar\rho>0$ such that $\delta\in C^2(\Omega_{\bar\rho})$.
If $\bar\rho$ is sufficiently small one can use $(\sigma,\delta)$
as new coordinates in $\Omega_{\rho}$, for all $\rho\in(0,\bar\rho)$.
In these coordinates the Laplacian becomes
$$\Delta_x=\frac{\partial^2}{\partial\delta^2}-
(N-1)\H\frac{\partial}{\partial\delta}+\Delta_{\sigma},$$
where $\Delta_{\sigma}$ is the Laplace-Beltrami operator on $\partial\Omega$
and $\H=\H(\cdot,\delta)$ is the mean curvature of $\Gamma_\delta$
(see \cite{BM-CV} for a detailed discussion).

Let $\eta=\eta(\delta)$ be a positive super-harmonic of $\L_\mu$ in $\Omega_\rho$,
as constructed in Lemma \ref{l-A}.
Set
$\overline\H:=(N-1)\,{\sup}_{\Omega_\rho}\H$.
Consider the initial value problem
\begin{equation}\label{e-ODE}
-\ddot{v}-\left(2\frac{\dot\eta}{\eta}-\overline\H\right)\dot{v}+\frac{\eta^{p-1}}{\delta^s}v^p=0
\qquad v(\rho)=0,\quad\dot v(\rho)=-\kappa,
\end{equation}
where $\kappa>0$.
Let $v_\kappa=v_\kappa(r)$ be the maximal left solution of \eqref{e-ODE} defined
on the maximal left interval of existence $(R_\kappa,\rho)$ in the region $\{(r,v)\in(0,\rho)\times\R\}$
(cf. \cite[pp.\,10-12 and 24-36]{Protter}).

Observe that $v_\kappa>0$ and $\dot v_\kappa<0$ for all $r\in(R_\kappa,\rho)$.
Indeed, if $r_0=\max\{r\in(R_\kappa,\rho):\dot v_\kappa(r)=0\}$ then $\ddot v_\kappa(r_0)>0$.
As $\kappa>0$, we conclude that $\{r\in(R_\kappa,\rho):\dot v_\kappa(r)=0\}=\emptyset$
and $v_\kappa(r)$ is {\sl strictly decreasing} on any interval.
In particular,
\begin{equation}\label{const-1}
\liminf_{r\to 0}v_\kappa(r)>0.
\end{equation}
An important consequence of the monotonicity of the solutions $v_\kappa$
is that they can be used to construct super-solutions of \eqref{*}.
\begin{lemma}\label{1d}
Let $\eta=\eta(\delta)$ be a positive super-harmonic of $\L_\mu$ in $\Omega_\rho$,
and $v_\kappa:(R_\kappa,\rho)\to\R$ be the maximal left solution of \eqref{e-ODE}.
Then $z_\kappa(\delta):=v_\kappa(\delta)\eta(\delta)$ is a super-solution to \eqref{*} in $\Omega_{R_\kappa,\rho}$.
\end{lemma}

\begin{proof}
A direct computation (cf. \cite[p.8]{Protter}) using the monotonicity of $v_\kappa$ shows that
\beq
\L_\mu z_{\kappa}&=& -\ddot{v}_\kappa\eta-v_\kappa\ddot{\eta}-2\dot{v}_\kappa\dot{\eta}+
(N-1)\H\left(\dot{v}_\kappa\eta+v_\kappa\dot{\eta}\right)-\frac{\mu}{\delta^2}v_\kappa\eta\\
&\ge&\left(-\ddot{v}_\kappa-\left(2\frac{\dot\eta}{\eta}-\overline\H\right)\dot{v}_\kappa\right)\eta
+\left(\L_\mu\eta\right)v_\kappa\\
&\ge&-\frac{\eta^p}{\delta^s}v_\kappa^p\quad\text{in }\:\Omega_{R_\kappa,\rho},
\eeq
as required.
\end{proof}

Our analysis of \eqref{e-ODE} is be based on the following well known ODE comparison lemma,
which we present here for reader's convenience.
\begin{lemma}\label{ODE-comp}
Assume that $u>0$ and $v>0$ satisfy differential inequalities
$$-\ddot{u}-a(r)\dot{u}+b(r)u^p\ge 0,\qquad
-\ddot{v}-a(r)\dot{v}+b(r)v^p\le 0
\qquad (r\in(R,\rho)),
$$
where $a,b\in C(R,\rho)$, $b\ge 0$ and $p>1$.
Then
\begin{itemize}
\item[(i)] {\em (IVP)-case}:
$u(\rho)\le v(\rho)$ and $\dot{u}(\rho)>\dot{v}(\rho)$ imply $u(r)<v(r)$ for all $r\in(R,\rho)$;
\item[(ii)] {\em (BVP)-case}:
$u(\rho)>v(\rho)$ and $u(R)>v(R)$ imply $u(r)>v(r)$ for all $r\in(R,\rho)$.
\end{itemize}
\end{lemma}

\begin{proof}
Part $(i)$ could be proved similarly to \cite[pp.\,26]{Protter}.
Part $(ii)$ can be established following the arguments in the proof of Lemma \ref{comparison}.
\end{proof}

Lemma \ref{KO-impr} follows via Lemma \ref{1d} from the following.

\begin{lemma}\label{l-ODE}
{\sc (ODE Lemma)}
Let $\mu\le 1/4$, $p>1$ and $s\ge\beta_-(p-1)+2$. If $v_\kappa$ is the
maximal left solution of \eqref{e-ODE} on the maximal existence interval
$(R_\kappa,\rho)$ then
\begin{itemize}
\item[$(i)$]
$R_\kappa>0$ and $v_\kappa(r)\to+\infty$ as $r\searrow R_\kappa$;
\item[$(ii)$]
$R_\kappa\to 0$ as $\kappa\to 0$;
\item[$(iii)$]
for any $r_\ast\in(0,\rho)$ one has
$\sup_{[r_\ast,\rho]}v_\kappa\to 0$ as $\kappa\to 0$.
\end{itemize}
\end{lemma}

\begin{proof}
To prove the lemma, one only has to show that $R_\kappa>0$.
As $v_\kappa(r)$ is decreasing in $r$ this obviously implies $v_\kappa(r)\to+\infty$ as $r\searrow R_\kappa$.

Indeed, assume that $(i)$ holds. Let $0<\kappa_1<\kappa_2$.
Then $v_{\kappa_1}<v_{\kappa_2}$ for all $r\in(R_{\kappa_2},\rho)$ by Lemma \ref{ODE-comp} $(ii)$.
In particular, this implies that $R_{\kappa_1}\le R_{\kappa_2}$.

Fix $r_\ast\in(0,\rho)$. For $\epsilon>0$, let $v^{(\epsilon)}$ be the unique solution
of the boundary value problem
\begin{equation}\label{e-BDE}
-\ddot v-\left(2\frac{\dot\eta}{\eta}-\overline\H\right)\dot v+\frac{\eta^{p-1}}{r^s}v^p=0,
\qquad v(r_\ast)=\epsilon,\quad v(\rho)=0.
\end{equation}
Set $-\kappa(\epsilon)=\dot v^{(\epsilon)}(\rho)$.
Thus $v_{\kappa(\epsilon)}=v^{(\epsilon)}$ for $r\in(r_\ast,\rho)$
in view of the uniqueness of solution for both \eqref{e-ODE} and \eqref{e-BDE}.
Moreover, $v_{\kappa(\epsilon)}\le\epsilon$ for $r\in(r_\ast,\rho)$
as $v_{\kappa(\epsilon)}$ is decreasing and
$\kappa(\epsilon)$ is strictly decreasing in view of the
BVP-comparison principle of Lemma \ref{ODE-comp} for equation \eqref{e-ODE}.
This proves $(ii)$ and $(iii)$.
\medskip

Now we are going to show that $R_\kappa>0$ for all $\kappa>0$. To do
this, we shall consider separately the cases $\mu<1/4$ and $\mu=1/4$,
with different choices of the super-harmonics $\eta$.

\paragraph{{\it Case $\mu<1/4$.}}
Here we choose a super-harmonic $\eta(r):=r^{\beta_-}(1+r^{\epsilon})$ and $\epsilon\in(0,1)$
(see Lemma \ref{l-A} (i)).
Then \eqref{e-ODE} can be written as
\begin{equation}\label{e1-ODE}
-\ddot v-\frac{2\beta_-}{r}\left(1+O(r^{\epsilon})\right)\dot v+
r^{(p-1)\beta_--s}(1+r^{\epsilon})^{p-1}v^p=0.
\end{equation}
Assume that $R_\kappa=0$ for some $\kappa>0$.
A direct computation (similar to the one in Proposition~\ref{supersol}) shows that
for a sufficiently large constant $\gamma>0$ and all $R\in(0,\rho)$
$$\ol v_R=\gamma r^{\frac{s}{p-1}-\beta_-}(r-R)^{-\frac{2}{p-1}}$$
is a super-solution to \eqref{e-ODE} in $(R,\rho)$, with $\gamma$ independent of $R$.
By Lemma \ref{ODE-comp} $(i)$ we conclude that
\begin{equation}\label{KO-ODE}
v_\kappa\le \gamma r^{\frac{s-2}{p-1}-\beta_-}\quad\text{in }\:(0,\rho).
\end{equation}
In the {\sl subcritical case} $s>\beta_-(p-1)+2$ this bound contradicts to \eqref{const-1},
so we conclude that $R_\kappa>0$.

In the {\sl critical case} $s=\beta_-(p-1)+2$,
linearizing \eqref{e-ODE} on $v_\kappa$ and taking into account \eqref{const-1}
we conclude that $v_\kappa$ is a sub-harmonic to the equation
\begin{equation}\label{e1-0}
-\ddot v-\frac{2\beta_-}{r}\left(1+O(r^{\epsilon})\right)\dot v+
\frac{C(1+r^{\epsilon})^{p-1}}{r^2}v=0\quad\text{in }\:(0,\rho/2),
\end{equation}
where $C:=\inf_{(0,\rho/2)}v_\kappa^{p-1}>0$. Let $\alpha_-<\alpha_+$ be the roots of the quadratic equation
$$\alpha(\alpha+2\beta_- -1)=C.$$
Note that $\alpha_-<0$ 
as $\beta_-<1/2$,
and choose $\alpha_-^\prime\in(\alpha_-,0)$.
A direct computation shows that for some $\rho_1\in(0,\rho/2)$ the function $\ol{h}:=Ar^{\alpha_-^\prime}$
is a super-harmonic to \eqref{e1-0} on $(0,\rho_1)$.
Choose $A>0$ in such a way that $\ol{h}(\rho_1)<v_\kappa(\rho_1)$ and $\dot{\ol{h}}(\rho_1)>\dot v_\kappa(\rho_1)$.
Then
$$v_\kappa\ge\ol{h}$$
by Lemma \ref{ODE-comp} $(ii)$.
But this contradicts to \eqref{KO-ODE}, and we conclude that $R_\kappa>0$.

\paragraph{{\it Case $\mu=1/4$.}}
Choose a super-harmonic $\eta(r):=r^{1/2}\left(1-(\log(1/r)^{-\epsilon}\right)$
and $\epsilon\in(0,1)$ as in Lemma \ref{l-A} (i).
Then \eqref{e-ODE} can be written as
\begin{equation}\label{e2-ODE}
-\ddot v-\frac{1}{r}\left(1+O(\log(1/r)^{-\epsilon-1})\right)\dot v+
r^{(p-1)/2-s}(1-\log(1/r)^{-\epsilon})^{p-1}v^p=0.
\end{equation}
Assume that $R_\kappa=0$ for some $\kappa>0$.
A direct computation shows that for a sufficiently large constant $\gamma>0$ and all $R\in(0,\rho)$
$$\ol v_R=\gamma r^{\frac{s}{p-1}-\frac{1}{2}}(r-R)^{-\frac{2}{p-1}}$$
is a super-solution to \eqref{e2-ODE} in $(R,\rho)$, with $\gamma$ independent of $R$.
As in \eqref{KO-ODE}, we obtain
\begin{equation}\label{KO-ODE-2}
v_\kappa\le \gamma r^{\frac{s-2}{p-1}-\frac{1}{2}}\quad\text{in }\:(R,\rho),
\end{equation}
for all small $R>0$.
In the {\sl subcritical case} $s>\frac{p+3}{2}$ this bound contradicts to \eqref{const-1},
so we conclude that $R_\kappa>0$.

In the {\sl critical case} $s=\frac{p+3}{2}$,
we simply observe that $v_\kappa$ is a sub-harmonic to the homogeneous equation
\begin{equation}\label{e1-1}
-\ddot v-\frac{1}{r}\left(1+O(\log(1/r)^{-\epsilon-1})\right)\dot v=0\quad\text{in }\:(0,\rho).
\end{equation}
On the other hand, a direct computation shows that
the function $\ol{h}:=A\log^{1/2}(1/r)$
is a super-harmonic to \eqref{e1-0} on $(0,\rho_1)$, for some $\rho_1\in(0,\rho)$.
Choose $A>0$ in such a way that $\ol{h}(\rho_1)<v_\kappa(\rho_1)$ and
$\dot{\ol{h}}(\rho_1)>\dot v_\kappa(\rho_1)$. Then
$$v_\kappa\ge\ol{h}$$
by Lemma \ref{ODE-comp} $(ii)$.
But this contradicts to \eqref{KO-ODE-2}, and we conclude that $R_\kappa>0$.
\end{proof}

\subsubsection{ Existence}

To prove the existence part of Theorem \ref{t-Main},
we first establish the existence of a solution between ordered sub- and super-solutions.

\begin{lemma}\label{l-subsup}
Let $\mu\le 1/4$ and $p>1$.
Assume that \eqref{*} admits a sub-solution $\usub$ and a super-solution $\usup$ in $\Omega$
so that $0\le\usub\le\usup$ in $\Omega$. Then \eqref{*} has a solution $U$ in $\Omega$
such that $\usub\le U\le\usup$ in $\Omega$.
\end{lemma}

\begin{proof}
For small $\eps>0$, let $U_\eps$ be a positive solution of
$$
\L_\mu U_\eps + \frac{U_\eps^p}{\delta^s} = 0\quad\mbox{in }\:D_\eps,
\qquad U_\eps=\usub\quad\mbox{on }\:\partial D_\eps.
$$
Such a solution is obtained, e.g., by minimization of the convex, coercive functional
$$
\int_{D_\eps} |\nabla U|^2-\frac{\mu}{\delta^2}U^2+\frac{|U|^{p+1}}{(p+1)\delta^s}\,dx
$$
in $H^1(D_\eps)$ with $U=\usub$ on $\partial D_\eps$.
By applying the Comparison Principle of Lemma \ref{comparison} $(ii)$ we obtain
$\usub\le U_\eps\le\usup$ on $D_{\eps}$.
Applying interior regularity together with the usual diagonalization argument we conclude that
$U=\lim_{\eps\to 0} U_\eps$ is the required solution of \eqref{*} in $\Omega$.
\end{proof}

Now, we prove the existence of $XXL$-solution in all of $\Omega$.

\begin{lemma}\label{l-XXL}
Let $\mu\le 1/4$, $p>1$ and $s<\beta_-(p-1)+2$.
Then \eqref{*} admits an $XXL$-solution in $\Omega$.
\end{lemma}

\begin{proof}
Let $\mu\le1/4$. Set
$$\usub_\rho:=\gamma\left(\delta^\frac{2-s}{1-p}-\kappa\delta^{1/2}\log(1/\delta)^{1/2}\right),$$
where $\kappa>0$ is chosen in such a way that $\usub_\rho(\rho)=0$. For some $\rho>0$ and sufficiently small $\gamma>0$,
the function $\usub_\rho$ is a sub-solution to \eqref{*} in
$\Omega_\rho$, cf. Table \ref{ko_tab} and the fact that
$\delta^{1/2}\log(1/\delta)^{1/2}$ is a local super-harmonic to $\L_{1/4}$
and hence a local super-harmonic to $\L_\mu$ for all $\mu\leq 1/4$
cf. Lemma~\ref{l-A}(ii). Let $\usub$ denote the function $\usub_\rho$, extended by zero to $D_\rho$.
Thus $\usub\ge 0$ is a sub-solution to \eqref{*} in the entire domain $\Omega$.

Set $\usup:=\gamma_\ast d^\frac{2-s}{1-p}$, where $d$ is the
Whitney-distance.
Note that $\usub\le\usup$ in $\Omega$ in view of the Keller--Osserman bound of Proposition~\ref{ko}.
Moreover, $\usup$ is a super-solution to \eqref{*} in $\Omega$, according to Proposition~\ref{supersol}.
By Lemma~\ref{l-subsup} we conclude that \eqref{*} admits a solution $U$ in $\Omega$
so that $\usub\le U\le\usup$ in $\Omega$, which is the required $XXL$-solution.
\end{proof}

\begin{remark}
The constructed $XXL$-solution $U$ satisfies, for some $\gamma>0$,
\begin{equation}\label{XXL-bound}
\gamma\le\liminf_{x\to\partial\Omega}\frac{U}{\delta^\frac{2-s}{1-p}}\le
\limsup_{x\to\partial\Omega}\frac{U}{\delta^\frac{2-s}{1-p}}\le\gamma^{-1}.
\end{equation}
\end{remark}

Next, we prove the existence of an $ML$-solution in all of $\Omega$.

\begin{lemma}\label{l-ML}
Let $\mu\le 1/4$, $p>1$ and $s<\beta_-(p-1)+2$.
Then \eqref{*} admits an $ML$-solution $U$ in $\Omega$.
\end{lemma}

\begin{proof}
We consider separately the cases $\mu<1/4$ and $\mu=1/4$.

\noindent
{\it Case $\mu<1/4$.}
Let $\alpha\in(\beta_-,\min\{\beta_-p+2-s,\beta_-+1,\beta_+\})$ and $\kappa>0$.
Set
$$
\usub_\rho:=\delta^{\beta_-}-\kappa\delta^{\alpha},
$$
where $\kappa>0$ is chosen in such a way that $\usub_\rho(\rho)=0$.
A direct computation shows that for a sufficiently small $\rho>0$,
$$\L_\mu\usub_\rho+\delta^{-s}\usub_\rho^p\leq
-\kappa(\alpha(1-\alpha)-\mu)\delta^{\alpha-2}(1+o(1))
+\delta^{-s}\left(\delta^{\beta_-}-\kappa\delta^{\alpha}\right)^p\le 0\quad\text{in }\:\Omega_{\rho},$$
that is $\usub_\rho$ is a sub-solution of \eqref{*} in $\Omega_\rho$.
Let $\usub$ denote the function $\usub_\rho$, extended by zero to $D_\rho$.
Hence $\usub\ge 0$ is a sub-solution to \eqref{*} in the entire domain $\Omega$.

Fix $\epsilon\in(0,\min\{1,\sqrt{1-4\mu}\})$. Then
$\Hsup:=\delta^{\beta_-}(1+\delta^\epsilon)$
is a large local super-harmonic of $\L_\mu$, as constructed in Lemma~\ref{l-A}.
We may assume that $\L_\mu\Hsup\ge 0$ in $\Omega_\rho$
(otherwise we adjust $\rho$ in the construction of $\usub$).
Let $R\in(0,\rho)$.
Let $\usup_R = \gamma_\ast d^{\frac{s}{p-1}}(d-R)^{-\frac{2}{p-1}}$, where $\gamma_\ast>0$
is chosen in such a way that $\usup_R$ is a super-solution to \eqref{*} in $D_{R}$,
see Lemma \ref{supersol}. Choose $\tau_\ast>1$ large enough, so that
$\tau_\ast\Hsup>\usup_R$ on $\Gamma_\rho$. Then
$$\usup:=\min\{\tau_\ast\Hsup,\usup_R\}$$
is a super-solution to \eqref{*} in the entire $\Omega$.

Note that $\usub\le\usup$ in $\Omega$, in view of the Comparison Principle of Lemma \ref{comparison} $(i)$.
By Lemma~\ref{l-subsup} we conclude that \eqref{*} has a solution $U$ in $\Omega$
so that $\usub\le U\le\usup$ in $\Omega$, which is the required $ML$-solution.

\noindent{\it Case $\mu=1/4$.}
Let $\alpha\in(0,1)$ and $\kappa>0$.
Set
$$\usub_{\rho}:=\delta^{1/2}\log\frac{1}{\delta}-
\kappa\delta^{1/2}\log^\alpha\frac{1}{\delta},$$
where $\kappa>0$ is chosen in such a way that $\usub_{\rho}(\rho)=0$.
A direct computation shows that for a sufficiently small $\rho>0$,
\begin{eqnarray*}
\L_{1/4}\usub_\rho+\delta^{-s}\usub_\rho^p &\leq&
-\kappa\alpha(1-\alpha)\delta^{-3/2}\left(\log^{\alpha-2}\frac{1}{\delta}\right)(1+o(1))\\
&&+\delta^{-s}\left(\delta^{1/2}\log\frac{1}{\delta}-\kappa\delta^{1/2}\log^\alpha\frac{1}{\delta}\right)^p\le 0
\end{eqnarray*}
in $\Omega_{\rho}$, that is $\usub_\rho$ is a sub-solution of
\eqref{*} in $\Omega_\rho$.

To construct a super-solution to \eqref{*}, fix $\epsilon\in(0,1)$ and set
$$\Hsup:=\delta^{1/2}\log\frac{1}{\delta}
\left(1+\log^{-\epsilon}\frac{1}{\delta}\right).$$
Thus $\Hsup$ is a large local super-harmonic of $\L_{1/4}$, see Lemma \ref{l-A}.
The rest of the proof is similar to the case $\mu<1/4$ above, so we omit it.
\end{proof}

\begin{remark}
The constructed $ML$-solution $U$ satisfies the bound
$$1\le\liminf_{x\to\partial\Omega}\frac{U}{\Hsup}\le\limsup_{x\to\partial\Omega}\frac{U}{\Hsup}\le\tau_\ast.$$
\end{remark}

\section{$S$--solutions and solutions for arbitrary $\mu>1/4$}\label{small_conclud}

It is easy to see that equation \eqref{*} admits local $S$-subsolutions for all $p>1$, $s\in\R$ and $\mu\le 1/4$.
Below we are going to show that the existence of global $S$-solutions is controlled by the global Hardy constant
$C_H(\Omega)$ rather then by relations between $p$, $s$ and $\mu$.

\begin{theorem}
Let $\mu\le C_H(\Omega)$, $p>1$ and $s\in\R$.
Then \eqref{*} has no nontrivial $S$--subsolution in $\Omega$.
\end{theorem}

\begin{proof}
Let $\usub\ge 0$ be a nontrivial $S$--subsolution of \eqref{*} in $\Omega$.
Set $h_\ast:=\delta^{1/2}\log^{1/2}(1/\delta)$.
Note that for all $\mu\le1/4$, $h_\ast$ is a local super-harmonic of $\L_\mu$ and
$$
\lim_{x\to\partial\Omega}\frac{\usub}{h_\ast}=0,
$$
cf. Theorem \ref{table}(i). For $\kappa>0$, consider the family $\vsub_\kappa:=(\usub-\kappa h_\ast)_+$.
Clearly, $\vsub_\kappa\in H^1_c(\Omega)$ and $\L_\mu\vsub_\kappa\leq 0$ in $\Omega$.
Testing this inequality with $\vsub_\kappa$ yields
$$
C_H(\Omega)\int_\Omega \frac{\vsub_\kappa^2}{\delta^2}\:dx
\leq\int_\Omega |\nabla\vsub_\kappa|^2\:dx
<\mu\int_\Omega \frac{\vsub_\kappa^2}{\delta^2}\:dx,
$$
which means $\vsub_\kappa=0$ in $\Omega$ for every $\kappa>0$. We conclude that $\usub=0$.
\end{proof}

The following lemma is crucial in our construction of global solutions for $\mu>C_H(\Omega)$.

\begin{lemma}\label{small}
Let $\mu>C_H(\Omega)$, $p>1$ and $s\in\R$.
Then there exists $\rho>0$ such that for every $\eps\in(0,\rho)$ equation \eqref{*} in $D_\eps$
admits a positive solution $u_\eps\in H^1_0(D_\eps)$.
Moreover, $u_\eps\in H^1_0(D_\eps)$ and $u_\eps(x)$ is monotone nondecreasing as $\eps\to 0$.
\end{lemma}

\begin{proof}
For a small $\eps>0$, consider the problem
\begin{equation}\label{De}
\L_\mu u_\eps + \frac{u_\eps^p}{\delta^s} = 0,
\quad u_\eps\in H^1_0(D_\eps),
\end{equation}
and the corresponding functional
$$
J_\eps(u)=\int_{D_\eps}\frac{1}{2}|\nabla u|^2-\frac{\mu}{2\delta^2}u^2+
\frac{u_+^{p+1}}{(p+1)\delta^s}\,dx
$$
in $H^1_0(D_\eps)$.
It is standard to see that $J_\eps$ is coercive and weakly lower semicontinuous on $H^1_0(D_\eps)$.
Moreover, minimizers of $J_\eps$ are nonnegative and solve \eqref{De}.

Let $u_\eps\ge 0$ be a minimizer of $J_\eps$.
From the definition of Hardy's constant $C_H(\Omega)$, it follows that if $\mu>C_H(\Omega)$ then
$u=0$ is not a local minimum of $J_\eps$ for $\eps>0$ sufficiently small.
Hence $u_\eps\gneq 0$ is the required solution of \eqref{De}.

Further, by applying the Comparison Principle of Lemma \ref{comparison} $(ii)$ we conclude
that $u_\eps(x)$ is monotone nondecreasing as $\eps\to 0$.
\end{proof}

\begin{theorem}
Let $\Omega$ be such that $C_H(\Omega)<1/4$.
Let $\mu\in(C_H(\Omega),1/4]$, $p>1$ and $s\in\R$.
Then equation \eqref{*} admits a positive $S$-solution in $\Omega$.
\end{theorem}

\begin{proof}
Let $\hsup>0$ be a super-harmonic in $\Omega_\rho$ for some $\rho>0$, as constructed in Lemma~\ref{l-A}.
For some fixed $R\in(0,\rho/2)$, let $\usup_R = \gamma_\ast d^{\frac{s}{p-1}}(d-R)^{-\frac{2}{p-1}}$
be a super-solution to \eqref{*} in $D_{R}$, as constructed in Lemma \ref{supersol}.
Choose $\tau_\ast\ge 1$ large enough, so that $\tau_\ast\hsup>\usup_R$ on $\Gamma_{\rho/2}$. Then
$$
\usup:=\min\{\tau_\ast\hsup,\usup_R\}
$$
is a super-solution to \eqref{*} in the entire $\Omega$.

Let $u_\eps$ be the monotone increasing family of solutions \eqref{*} in $D_\eps$, as constructed in Lemma \ref{small}.
By applying the Comparison Principle of Lemma \ref{comparison} $(ii)$ we obtain
$u_\eps\le\usup$ on $D_{\eps}$.
Applying the usual diagonalization argument we conclude that
$u:=\lim_{\eps\to 0} u_\eps$ is the required $S$-solution of \eqref{*} in $\Omega$.
\end{proof}

\begin{remark}
Observe that if $s<\beta_+(p-1)+2$ and $C_H(\Omega)<\mu\leq 1/4$ then
$$\lim_{x\to\partial\Omega}\frac{\delta^{\frac{s-2}{p-1}}}{\hsup}=0$$
for every positive local super-harmonic $\hsup$ of $\L_\mu$, see Theorem \ref{table-sup}.
By Lemma \ref{supersol} and the Comparison Principle of Lemma \ref{comparison} $(ii)$ we obtain
that if $s<\beta_+(p-1)+2$ then every $S$-subsolution $u$ of \eqref{*} satisfies an improved upper bound
$$u\le \gamma_\ast \delta^{\frac{s-2}{p-1}}\quad\text{in }\:\Omega,$$
which is stronger then the upper bound on $S$-subsolutions imposed by positive super-harmonics.
\end{remark}

Our classification of (sub)\,solutions to \eqref{*} is not applicable for $\mu>1/4$.
However, one can show that for all values of $\mu>1/4$,
equation \eqref{*} admits positive solutions which obey the Keller--Osserman bound.

\begin{theorem}\label{super-mu}
Let $\mu>1/4$ and $p>1$.
Then equation \eqref{*} admits a positive solution $u$ in $\Omega$
such that $u\le\gamma_\ast\delta^{\frac{2-s}{1-p}}$ in $\Omega$.
\end{theorem}

\begin{proof}
Let $u_\eps$ be the monotone increasing family of solutions \eqref{*} in $D_\eps$, as constructed in Lemma \ref{small}.
By Lemma \ref{supersol} and the Comparison Principle of Lemma \ref{comparison} $(ii)$ we obtain
$u_\eps\le \gamma_\ast d^{\frac{s-2}{p-1}}$ on $D_{\eps}$.
Applying the usual diagonalization argument we conclude that
$u:=\lim_{\eps\to 0} u_\eps$ is the required solution of \eqref{*} in $\Omega$.
\end{proof}

\section{Open problems}\label{open}

We finish our investigation with a list of open problems, which we
consider as interesting:

\begin{problem}
If we assume $p>1$, $\mu\leq 1/4$ and $s<\beta_-(p-1)+2$ then in
Theorem \ref{t-Main} we have proved the existence of an $XXL$-solution
with boundary behaviour given by \eqref{XXL-bound}.
What is the precise boundary behaviour of an $XXL$-solution? We conjecture that the correct asymptotic behaviour is
given by $\lim_{x\to\partial\Omega}
U(x)/\delta(x)^\frac{s-2}{p-1}=\const$, where the constant depends
only on $p$, $s$ and $\mu$. In the case $s=0$, $\mu=0$ this was proved in \cite{BM-JAM}, \cite{BM-CV}
and \cite{MV2}.
\end{problem}

\begin{problem}
In the case $p>1$, $\mu=s=0$ every $XL$-solution is automatically an $XXL$-solution
and moreover the $XXL$-solution is unique, see \cite{BM-JAM,MV2}. Are these
two statements true for every $s$, $\mu$ in the existence range?
\end{problem}

\begin{problem}
What is the asymptotic behavior near the boundary of the solutions,
constructed in Theorem \ref{super-mu} for arbitrary $\mu>1/4$ ?
\end{problem}

\begin{problem}
Is the existence and non-existence threshold phenomena similar to Theorem \ref{t-Main}
valid for some (or maybe all) $p<1$, or is there a natural reason,
why the result can only be true for $p > 1$ ?
\end{problem}


\subsection*{Acknowledgements.}
The authors wish to thank Vitali Liskevich for valuable discussions.
The work was supported by the Royal Society grant "Liouville theorems in nonlinear elliptic equations and systems".
Part of this research was done while V.M.
was visiting the Universities of Basel and Z\"urich, and C.B. and W.R. were visiting
the University of Bristol. The authors would like to thank
these institutions for their kind hospitality. The work of W.R. was supported by a grant from the Swiss National Science Foundation.

\end{document}